\documentclass{amsart}

\usepackage{amsmath,amssymb,amsthm,latexsym}

\newtheorem{theorem}{Theorem}
\newcommand{\bt}{\begin{theorem}}
\newcommand{\et}{\end{theorem}}
\newtheorem{lemma}{Lemma}
\newcommand{\bl}{\begin{lemma}}
\newcommand{\el}{\end{lemma}}
\newtheorem{corollary}{Corollary}
\newcommand{\bc}{\begin{corollary}}
\newcommand{\ec}{\end{corollary}}
\newcommand{\bconj}{\begin{conjecture}}
\newcommand{\econj}{\end{conjecture}}
\newtheorem{problem}{Problem}
\newcommand{\bprob}{\begin{problem}}
\newcommand{\eprob}{\end{problem}}
\newcommand{\beq}{\begin{equation}}
\newcommand{\eeq}{\end{equation}}
\newcommand{\benum}{\begin{enumerate}}
\newcommand{\eenum}{\end{enumerate}}
\newcommand{\N}{\ensuremath{ \mathbf N }}
\newcommand{\Z}{\ensuremath{\mathbf Z}}

\newcommand{\mci}{\ensuremath{ \mathcal I}}

\newcommand{\mcr}{\ensuremath{ \mathcal R}}

\newcommand{\mcx}{\ensuremath{ \mathcal X}}

\newcommand{\mba}{\ensuremath{ \mathbf a}}

\newcommand{\mbx}{\ensuremath{ \mathbf x}}

\newcommand{\bmat}{\left(\begin{matrix}}
\newcommand{\emat}{\end{matrix}\right)}
\newcommand{\bsmallmat}{\left(\begin{smallmatrix}}
\newcommand{\esmallmat}{\end{smallmatrix}\right)}
\DeclareMathOperator{\qand}{\quad\text{and}\quad}
\DeclareMathOperator{\qqand}{\qquad\text{and}\qquad}

\title{Explicit sumset sizes in additive number theory}
\author{Melvyn B.  Nathanson}
\address{Department of Mathematics\\Lehman College (CUNY)\\Bronx, NY 10468}
\email{melvyn.nathanson@lehman.cuny.edu}

\date{\today}

\thanks{Supported in part by  PSC-CUNY Research Award Program grant 66197-00 54.}

\begin{document}

\begin{abstract}
It is an open problem in additive number theory to compute 
and understand  the full range of sumset sizes of finite sets 
of integers, that is, the set 
$\mcr_{\Z}(h,k) = \{|hA|:A \subseteq \Z \text{ and } |A|=k\}$ 
for all integers $h \geq 3$ and $k \geq 3$.   
This paper constructs certain infinite families of finite sets 
of size $k$, computes their $h$-fold sumset sizes, and obtains explicit
finite arithmetic progressions of sumset sizes in $\mcr_{\Z}(h,k)$.  
\end{abstract}

\maketitle

\section{Sums of finite sets of integers}
An additive abelian semigroup is simply a nonempty set $G$ with a commutative 
and associative binary operation, written additively.   
For every nonempty subset $A$ of the semigroup $G$, 
the \emph{$h$-fold sumset of $A$} is the set of all sums of $h$ 
not necessarily distinct elements of $A$, that is, 
\[
hA = \left\{a_{j_1}+\cdots + a_{j_h}: \text{$a_{j_r} \in A$ for all $r= 1,\ldots, h$} \right\}.
\]  
We define $0A = \{0\}$. 
A core part of additive number theory is the study of sumsets of finite subsets 
of additive abelian semigroups.  
We define the \emph{sumset size set} 
\[
\mcr_G(h,k) = \left\{ |hA|: A \subseteq G \text{ and } |A| = k \right\}.
\] 
A basic problem is to understand this set.

We  consider  the set $\mcr_{\Z}(h,k)$ of sumset sizes of finite sets of integers.  
The dilation of a set $A$ by $\lambda$ is the set 
$ \lambda\ast A =  \{\lambda a  : a \in A\}$. 
Sets $A$ and $B$ are \emph{affinely equivalent} if there are numbers 
$\lambda \neq 0$ and $\mu$ such that
\[
B = \lambda\ast A + \mu = \{\lambda a + \mu: a \in A\}.
\]
If $A$ and $B$ are affinely equivalent, then 
\[
hB = h(\lambda\ast A + \mu) = \lambda \ast hA + h\mu
\]
and so $|hB| = |hA|$. 
Because sumset size is an affine invariant, we have $\mcr_{\Z}(h,k) = \mcr_{\N_0}(h,k)$.  
It is proved in~\cite{nath25cc} that $\mcr_{\Z^n}(h,k) = \mcr_{\Z}(h,k)$ 
for all positive integers $n$. 

There are simple lower and upper bounds for $\mcr_{\Z}(h,k)$.  We have 
\beq          \label{explicit::min}
\min \mcr_{\Z}(h,k) = hk-h+1
\eeq
and, if $ A \subseteq \Z \text{ and } |A| = k $, then $|hA| = hk-h+1$
 if and only if $A$ is an arithmetic progression of length $k$. 
Similarly, if $A$ is a $B_h$-set, that is, a set all of whose $h$-fold sums are distinct, 
then $|hA| = \binom{h+k-1}{h}$ and 
\beq          \label{explicit::max}
\max \mcr_{\Z}(h,k) = \binom{h+k-1}{h}.
\eeq
Beginning with the work of Grigori Freiman~\cite{frei73,frei87,nath96bb}, 
a large research literature has investigated 
finite sets whose sumsets are very small, that is, close to the minimum size.  
There is also a large research literature~\cite{obry04} on sets whose sumsets are 
close to the maximum size.  

What is surprising is the lack of attention to the \emph{full range of sumset sizes}  of finite sets of integers.  
Possibly, the only published statement related to this problem occurs in a 
1983 paper by  Erd\H os and  Szemer\' edi~\cite{erdo-szem83} about the number 
of sums and products.  
They wrote: 
\begin{quotation}
Let $2k-1 \leq t \leq \frac{k^2+k}{2}$.  It is easy to see that one can find a sequence 
of integers $a_1 < \ldots < a_k$ so that there should be exactly $t$ distinct integers 
in the sequence $a_i+a_j, 1 \leq i \leq j \leq k$.  
\end{quotation}
Theorem~\ref{explicit::theorem:calc-2-k} refines this assertion.  
For real numbers $u$ and $v$, define the \emph{integer interval} 
\[
[u,v] = \{ n \in \Z: u \leq n \leq v\}.
\]

\bt[Nathanson~\cite{nath25bb}]                 \label{explicit::theorem:calc-2-k}                
For all positive integers $k$, 
\[
\mcr_{\Z}(2,k) =  \left[ 2k-1, \binom{k+1}{2} \right].
\]
Moreover, for all $t \in \mcr_{\Z}(2,k)$, there exists a set $A \subseteq \left[0,2^k -1 \right]$ 
such that $|A| = k$ and $\left| 2A \right| = t$. 
\et 

Here are two important observations.  First, the set of sumset sizes $\mcr_{\Z}(2,k)$ is known 
exactly, and it is an integer interval:  There is no ``missing number''   
between $\min \mcr_{\Z}(2,k)$ 
and $\max \mcr_{\Z}(2,k)$.  Second, there is a finite, albeit exponential, upper bound 
on the amount of computation needed to find a set $A$ with $|A| = k$ and $|2A| = t$ 
for all $t \in \mcr_{\Z}(2,k)$.  
For all $h \geq 2$, the set $\mcr(h,k)$ is finite and so there exists an integer $N$ such that, 
for all $t \in \mcr(h,k)$, there is a set $A$ with $|A| = k$, $|hA| = t$, and 
$A \subseteq [0,N]$. Let $N(h,k)$ be the least such number.  
By Theorem~\ref{explicit::theorem:calc-2-k}, $N(2,k) < 2^k$. 
For $h \geq 3$, there is the following exponential upper bound. 

\bt[Nathanson~\cite{nath25cc}]   
For all $h \geq 3$ and $k \geq 3$,
\[
N(h,k) < 4(4h)^{k-1}.
\]
\et
It would be of interest to know if, for fixed $h$, there is a subexponential 
or even polynomial upper bound for $N(h,k)$.

We have 
\[
\mcr_{\Z}(h,1) = \{1\}   \qand  \mcr_{\Z}(h,2) = \{h+1\} \qquad \text{for all $h \geq 1$}
\]
and 
\[
 \mcr_{\Z}(1,k) = \{k\} \qquad \text{for all $k \geq 1$.}
\]
Theorem~\ref{explicit::theorem:calc-2-k} describes $ \mcr_{\Z}(2,k)$. 
Thus, the problem is to understand $\mcr_{\Z}(h,k) $ for $h \geq 3$ and $k \geq 3$.

Let $k = 3$.  
From~\eqref{explicit::min} and~\eqref{explicit::max} we have 
\[
 \mcr_{\Z}(h,3) \subseteq \left[ 2h+1, \binom{h+2}{2}\right].
\]
In particular, 
\[
\mcr_{\Z}(3,3) \subseteq \{7,8,9,10\}.
\]
We have  
\begin{align*}
3\{0,1,2\} = [0,6]  & \qqand |3\{0,1,2\} | =  7 \\
3\{0,1,3\} = [0,7] \cup \{9\} & \qqand |3\{0,1,3\} | =  9 \\
3\{0,1,4\} = [0,6] \cup \{8,9,12\}  & \qqand |3\{0,1,4\} | =  10
\end{align*}
and so   
\[
\{7,9,10\} \subseteq \mcr_{\Z}(3,3).
\]
Where is 8?  A computer search failed to find a set $A$ of integers with $|A| = 3$
and $|3A| = 8$.  Nathanson~\cite{nath25bb}  proved that $8 \notin \mcr_{\Z}(3,3)$, 
that is, 
\[    
\mcr_{\Z}(3,3) = \{7,9,10\}.
\]
Why is there no 8 in $\mcr_{\Z}(3,3)$?  There is a proof but not a reason.

More generally, we have the following ``missing number'' result.  

\bt[Nathanson~\cite{nath25bb}]                                                            \label{explicit::theorem:calc-h-k}     
For all $h \geq 3$ and $k \geq 3$, 
\[
hk-h+2 \notin \mcr_{\Z}(h,k) 
\]
and so the sumset size set $\mcr_{\Z}(h,k)$ is not an interval.  
\et

Theorems~\ref{explicit::theorem:calc-2-k}  and~\ref{explicit::theorem:calc-h-k}  
inspire this field of research.  
For $h \geq 3$ and $k \geq 3$, the sumset size set  $\mcr_{\Z}(h,k)$ is not 
the integer interval defined by its minimum and 
maximum values.  What is it?  
One can generate random sets $A$ of size $k$, 
compute their sumset size $|hA|$, and generate subsets of the set $\mcr_{\Z}(h,k)$.  
From these tables one can formulate conjectures.  
Observation of gaps in tables of the sumset sizes of random subsets suggested 
the following result, which was proved by Vincent Schinina. 
 
\bt[Schinina~\cite{schi25}]
For all $h \geq 3$ and $k \geq 3$, 
\[
[hk-h+2,hk-1] \cap  \mcr_{\Z}(h,k) = \emptyset
\]
and 
\[
hk \in \mcr_{\Z}(h,k). 
\]
\et

For sets of size 3 there is the following result. 

\bt[Nathanson~\cite{nath25bb}] 
For all positive integers $h$, 
\beq                  \label{explicit:h-3} 
\mcr_{\Z}(h,3) = \left\{  \binom{h+2}{2} - \binom{\ell}{2} :  \ell \in [1,h] \right\}. 
\eeq  
\et

Thus, the sumset sizes of 3-element sets are differences of triangular numbers. 
We know $\mcr_{\Z}(h,3)$, but the set $\mcr_{\Z}(h,4)$ 
of sumset sizes of 4-element sets is still a mystery.  

While computation of sumset sizes of random finite sets of integers generates 
elements of the sumset size sets $\mcr_{\Z}(h,k)$, 
it is useful and important to have explicit constructions of finite sets 
and explicit families of sumset sizes in $\mcr_{\Z}(h,k)$.  
In Theorems~\ref{explicit:theorem:intervals}  and~\ref{explicit:theorem:abc} 
and Corollaries~\ref{explicit:corollary:intervals}--\ref{explicit:corollary:AP},   
we construct  infinite families of such sets, and, 
in particular, interesting elements of $\mcr_{\Z}(h,4)$.

\section{Arithmetic progressions of intervals}
Let $\N_0 = \{0,1,2,3,\ldots\}$ be the set of nonnegative integers.  
For positive integers $h$ and $k$, let  
\[
\mcx_{h,k} = \left\{ (x_1,\ldots, x_k) \in \N_0^k: \sum_{j=1}^k x_j = h\right\}.
\]
The set $\mcx_{h,k}$ is invariant under permutations:  For all $\sigma \in S_k$, 
we have $(x_1,\ldots, x_k) \in \mcx_{h,k}$ if and only if 
$(x_{\sigma(1)},\ldots, x_{\sigma(k)}) \in \mcx_{h,k}$. 

Let $A = \{a_1,\ldots, a_k\} \subseteq \Z$ with $|A| = k$ 
and let $\mba = (a_1,\ldots, a_k) \in \Z^k$. 
For $\mbx = (x_1,\ldots, x_k) \in \mcx_{h,k}$, we define 
\[
\mbx \cdot \mba = \left(x_1,\ldots, x_k)  \cdot (a_1,\ldots, a_k \right) 
= \sum_{j=1}^k x_ja_j.
\]
Then $\mbx \cdot \mba \in hA$ and 
\[
hA = \left\{ \mbx \cdot \mba : \mbx \in \mcx_{h,k} \right\}. 
\]
The vector \mba\ depends on the ordering of the elements of the set $A$,
but, because $\mcx_{h,k}$ is $S_k$-invariant, 
the sumset $hA$  is independent of the ordering of $A$. 

It is straightforward to check that,
for all  positive integers $h$ and $k$, 
\beq             \label{explicit:Xhk}
\mcx_{h,k} = \bigcup_{x_k=0}^{h} 
\left\{ (x_1,\ldots, x_{k-1}, x_k): (x_1,\ldots, x_{k-1}) \in \mcx_{h - x_k,k-1} \right\}.
\eeq

The following terminology is useful.  
Let $u_1 \leq u_2$.  We say that there is a \emph{gap} between integer intervals   
$[u_1,v_1]$ and $[u_2,v_2]$  if there is an 
integer $n$ such that 
\[
v_1 < n < u_2.  
\]
The integer intervals   $[u_1,v_1]$ and $[u_2,v_2]$ \emph{have no gap} if 
$u_2 \leq v_1+1$ or, equivalently if $[u_1,v_1] \cup [u_2,v_2] $ is an integer interval.

\bl                    \label{explicit:lemma:Ihk}
For all positive integers $h$ and $k$, 
\[
\mci_{h,k} =  \left\{  \sum_{j=2}^{k}     (j-1) x_j :  (x_1,\ldots, x_{k}) \in \mcx_{h,k} \right\} 
 = [0, (k-1) h ].
 \]
\el

\begin{proof} 
The proof is by induction on $k$. 
For $k=1$, we have $\mcx_{h,1} = \{(h)\}$ and $\mci_{h,1} = \{0\} = [0,0]$.  
For $k =2$, we have 
\begin{align*} 
 \mcx_{h,2}  
 &= \left\{ (x_1,x_2) \in \N_0^2: x_1+x_2= h\right\} \\ 
&  = \left\{ (h-x_2,x_2) \in \N_0^2: x_2 \in [0,h] \right\}
\end{align*} 
and so 
\[
\mci_{h,2} =  \left\{    x_2   :  (x_1,x_2) \in \mcx_{h,2} \right\}  = [0,h].
 \]
 
 Let $k \geq 3$ and assume that $\mci_{h,k-1} = [0, (k-2) h ]$ 
 for all $h \geq 1$. 
 
If $x_k \in [0,h-1]$, then 
\[
(k-1)(x_k+1) \leq (k-2)h+x_k+1.  
\]
It follows that there is no gap between the  integer intervals 
\[
 \left[ (k-1)x_k, (k-2)h +  x_k \right]
\]
and 
\[
 \left[ (k-1)(x_k +1), (k-2)h +  x_k + 1 \right]
\]
and so their union is the integer interval 
\begin{align*}
\left[ (k-1)x_k,  (k-2)h +  x_k + 1 \right]. 
\end{align*}
Applying relation~\eqref{explicit:Xhk}, we obtain 
\begin{align*}
\mci_{h,k}  
& =  \left\{  \sum_{j=2}^{k}     (j-1) x_j :  (x_1,\ldots, x_k) \in \mcx_{h,k} \right\} \\
& =  \left\{ (k-1)x_k+ \sum_{j=2}^{k-1} (j-1) x_j  : (x_1,\ldots, x_k) \in \mcx_{h,k} \right\} \\
& =  \bigcup_{x_k=0}^h \left\{ (k-1)x_k+ \left\{  \sum_{j=2}^{k-1}     (j-1) x_j : 
(x_1,\ldots, x_{k-1}) \in \mcx_{h - x_k,k-1} \right\}  \right\} \\
& =  \bigcup_{x_k=0}^h \left\{ (k-1)x_k+ \mci_{h - x_k,k-1} \right\} \\ 
& =  \bigcup_{x_k=0}^h \left\{ (k-1)x_k+ [0, (k-2) (h - x_k) ]\right\} \\ 
& =  \bigcup_{x_k=0}^h \left[ (k-1)x_k,  (k-2)h + x_k \right]. 
\end{align*}
Because there is no gap between consecutive pairs of these $h+1$  intervals, 
we obtain 
\[
\mci_{h,k} = [0, (k-1)h].
\]
This completes the proof. 
\end{proof}

\bt                \label{explicit:theorem:intervals} 
Let $k$, $a$, $b$, and $\ell$ be positive integers with $k=a\ell$ and $a \leq b$.  
Let 
\[
P = b\ast [0,\ell -1] = \{0,b,2b,\ldots, (\ell-1)b \} 
\] 
be the  $\ell$-term arithmetic progression  
with difference $b$ and smallest element 0. 
Let 
$A$ be the $\ell$-term arithmetic progression of translates of the interval $[0, a-1 ]$: 
\begin{align*}
A & = P+[0,a-1] 
= \bigcup_{j=1}^{\ell} \left(  (j-1)b + [0,a-1] \right).
\end{align*}
Then $|A| = k$.  

For every positive integer $h$, let 
\[
Q = b\ast [0, h(\ell -1)]   = \{0,b,2b,\ldots, h(\ell-1) b \} 
\] 
be the $(h(\ell-1) +1)$-term arithmetic progression  
with difference $b$ and smallest element 0. 
The sumset $hA$ is an 
$(h(\ell-1) +1)$-term arithmetic progression of translates of the interval $[0, h(a-1)]$: 
\beq           \label{explicit:hAstructure}
hA = Q +  \left[  0, h (a-1)   \right]  
\eeq
and 
\beq           \label{explicit:hAsize}
|hA| = \begin{cases}
(a+ b(\ell-1) - 1)h+1 & \text{if $a \leq b \leq (a-1)h+1$} \\
(a-1)(\ell-1)h^2 + (a+\ell-2)h +1 
& \text{if $b \geq h(a-1) +1$.}  
\end{cases}
\eeq
\et

\begin{proof}
We have 
\begin{align*}
A 
& = P+[0,a-1]  \\
& = \bigcup_{j=1}^{\ell}  [ (j-1)b, \ a -1 + (j-1)b ] \\ 
& = \bigcup_{j=1}^{\ell } L_j
\end{align*}
where $L_j$ is the integer interval 
\[
L_j =  [ (j-1)b, \  a  -1 + (j-1)b ] 
\]
and $|L_j|=a$.  The inequality $a \leq b$ implies $a -1 +(j-1)b < jb$ for all $j \in [1,\ell]$ 
and so the integer intervals $L_j$ 
are pairwise disjoint  and $|A| = \sum_{j=1}^{\ell} |L_j| = a\ell = k$.  

If $ (x_1,\ldots, x_{\ell}) \in \mcx_{h,\ell}$, then $\sum_{j=1}^k x_j = h$.  
Applying Lemma~\ref{explicit:lemma:Ihk} with $k=\ell$, we have 
\begin{align*}
hA & = h\left(\bigcup_{j=1}^{\ell } L_j \right)  
 = \bigcup_{ \mbx  = (x_1,\ldots, x_{\ell}) \in \mcx_{h,\ell} } 
  \left(x_1 L_1 + \cdots + x_{\ell}  L_{\ell} \right) \\
& = \bigcup_{ \mbx  \in \mcx_{h,\ell} } \sum_{j=1}^{\ell}  
 [ (j-1) x_j  b,   (a-1) x_j +  (j-1)  x_j b ] \\
& = \bigcup_{\mbx \in \mcx_{h,\ell}}  \left[  \sum_{j=1}^{\ell }   (j-1)  x_j b, 
  \sum_{j=1}^{\ell } (a-1)  x_j  +  \sum_{j=1}^{\ell }    (j-1) x_j b  \right] \\
& = \bigcup_{\mbx \in \mcx_{h,\ell}}  \left[  \sum_{j=2}^{\ell }  (j-1)  x_j b,  h(a-1)   
+ \sum_{j=2}^{\ell }      (j-1)  x_j b  \right] \\
& = \bigcup_{\mbx \in \mcx_{h,\ell}}   \left( b \sum_{j=2}^{\ell } (j-1)  x_j   
+  \left[  0,  h(a-1)  \right]  \right)  \\
& = \left\{ b  \sum_{j=2}^{\ell }    (j-1)  x_j  : \mbx \in \mcx_{h,\ell}\right\}  
+  \left[  0, h (a-1)  \right] \\
& = b\ast  \mci_{h,\ell}  +  \left[  0, h (a-1)   \right] \\
& = b\ast  [0, (\ell -1) h ]  +  \left[  0, h (a-1)   \right] \\
& = Q +  \left[  0, h (a-1)   \right]. 
\end{align*}
This proves~\eqref{explicit:hAstructure}. 

To obtain the sumset size formula~\eqref{explicit:hAsize}, we write 
$hA$ as a union of intervals:
\begin{align*}
hA 
&= \bigcup_{j=0}^{ (\ell -1) h} \left[ bj, bj + h (a-1)   \right]. 
\end{align*}
If $b \geq h(a-1)+1$, then the $(\ell -1) h +1$ intervals are pairwise disjoint and 
\[
|hA| = ((\ell -1) h +1)( (a-1)h+1).
\] 
If $a \leq b \leq h(a-1)+1$, then there are no gaps between successive intervals 
and so 
\[
hA = [0, b(\ell -1) h +(a-1)h] 
\] 
and 
\[
|hA| = (a+ b(\ell -1)   -1)h +1.
\] 
This completes the proof.  
\end{proof}

The following results are immediate consequences of Theorem~\ref{explicit:theorem:intervals}.

\bc                            \label{explicit:corollary:intervals} 
Let $h$ and $k$ be positive integers
If  $a$   and $\ell$ are positive integers such that $k=a\ell$, then 
the sumset size set $\mcr_{\Z}(h,k)$ contains the arithmetic progression 
\[
\left\{ (\ell -1) hb + (a -1)h +1 :b \in [a, (a-1)h+1] \right\}.
\]
\ec

\bc                        \label{explicit:corollary:Rh4} 
For every positive integer $h$, the sumset size set  $\mcr_{\Z}(h,4)$ 
contains the $h$-term arithmetic progression 
$\{bh+1:b\in [3,h+2] \}$.
\ec

\bc
If $k = a^2$, then $\left((a-1)h+1\right)^2 \in \mcr_{\Z}(h,k)$. 
\ec

\section{Sums of  intervals of different  lengths}
Consider a finite set that is the union of two intervals of different lengths.  
Let $a$, $b$, and $c$ be nonnegative  integers with $a < b$ and $a \neq c$, 
and let $A = [0,a] \cup [b,b+c]$.  
The set $A$ is affinely equivalent to the set 
\begin{align*}
A' & =   (-1)\ast A + b+c  \\
& = [0,c]  \cup [ b+c-a,  b+c] \\
& = [0,a'] \cup [b', b'+c']
\end{align*} 
with $a'=c <  b+c-a = b'$, and $c' = a$.  
If $a < c$, then  $a' > c'$.  
Moreover, $|A| = |A'|$ and  $|hA| = |hA'|$.     
Thus, it suffices to consider only the case $a>c$. 

 Note that the case $a=c$ (that is, $A = [0,a] \cup [b,b+a] = \{0,b\}+[0,a]$)   
 is a special case of Theorem~\ref{explicit:theorem:intervals}.

The \emph{integer part} (also called the \emph{floor}) of the real number $w$,  
denoted $[w]$, is the unique integer $n$ such that $n \leq w < n+1$. 
There should be no notational confusion between $[u,v]$ and $[w]$.

\bt               \label{explicit:theorem:abc}
Let $a$, $b$, and $c$  be integers with $0 \leq c < a < b$ and let  
\[
A= [0,a] \cup [b,b+c].  
\]
Let $h \geq 2$ and 
\beq                                                                  \label{explicit:i0}
i_0 = \left[\frac{ha-b}{a-c}\right]. 
\eeq
If $b> ha$, then 
\beq                                                                  \label{explicit:c>ha}
|hA| =  (h+1)\left( 1 + \frac{h(a+c)}{2}  \right).
\eeq
If $a < b \leq ha$, then 
\begin{align}                 \label{explicit:c<ha,i0<}
|hA| & = (i_0+1) b  +  (h-i_0)( ha + 1) - \frac{(h+i_0+1)(h-i_0)(a-c)}{2}.        
\end{align}
\et

\begin{proof}
We have 
\begin{align*}
hA & = \bigcup_{i=0}^h \left(  (h-i)  [0,a] + i[b,b+c]   \right) \\
& = \bigcup_{i=0}^h [ib, (h-i)a+ i(b+c)]  \\
& = \bigcup_{i=0}^h L_i
\end{align*}
where $L_i$ is the integer interval 
\begin{align*}
L_i & =  [ib, (h-i)a+  i(b+c) ] \\
& = [ib, ha +i(b-a+c)]
\end{align*} 
and 
\[
\left|  L_i \right| =  ha + 1 - i (a-c).
\]
Because $b \geq 1$, the lower bounds of the intervals $L_0, L_1,\ldots, L_h$ 
are strictly increasing.  
Let $i \in [0,h-1]$.  Because $a > c$,  
the intervals $L_i$ and $L_{i+1}$  are disjoint if and only if 
\[
 ha - i (a-c) + ib < (i+1)b 
\]
if and only if 
\beq          \label{explicit:i0}
\frac{ha-b}{a-c} < i \leq h-1. 
\eeq
Let
\[
i_0 = \left[ \frac{ha-b}{a-c}  \right].
\]
Thus, for $i \in [0,h-1]$, 
the intervals $L_i$ and $L_{i+1}$  are disjoint if and only if 
\[
i \in [i_0+1,h-1]. 
\]
Thus, the $h+1$ intervals $L_i$ are pairwise disjoint if and only if $i_0 \leq -1$ 
or, equivalently, if and only if $b > ha$. 
In this case, 
\begin{align*} 
|hA| 
& = \left| \bigcup_{i=0}^h L_i \right| = \sum_{i=0}^h  \left|  L_i \right| \\
& = \sum_{i=0}^h  \left(   ha  + 1 - i (a-c) \right) \\ 
& = (h+1) \left( ha + 1  \right) - \frac{h(h+1)(a-c)}{2} \\
& = (h+1)\left( 1 + \frac{h(a+c)}{2}  \right). 
\end{align*}

If $b \leq ha$ or, equivalently, if $i_0 \geq 0$, then the $h-i_0-1$ intervals 
$L_{i_0+2},L_{i_0+3}, \ldots, L_h$ are pairwise disjoint 
and $L_i \cap L_{i+1} \neq \emptyset$ for $i \in [0,i_0]$. 
It follows that 
\[
L_0^* = \bigcup_{i=0}^{i_0+1} L_i = [0, ha + (i_0+1)(b-a+c)]  
\] 
and  
\[
 \left| L_0^* \right| =  ha + 1 + (i_0+1)(b-a+c). 
\] 
Moreover,
\begin{align*}
 \left|  \bigcup_{i=i_0+2}^h L_i  \right| 
 & =  \sum_{i=i_0+2}^h |L_i |  = \sum_{i=i_0+2}^h(   ha  + 1- i (a-c)) \\ 
& = (h-i_0-1)(ha + 1) - \frac{(h+i_0+2)(h-i_0-1)(a-c)}{2} .
\end{align*}
and 
\[
L_0^* \cap \left( \bigcup_{i=i_0+2}^h L_i  \right) =
\left( \bigcup_{i=0}^{i_0+1} L_i \right) \cap \left( \bigcup_{i=i_0+2}^h L_i  \right) 
 = \emptyset.  
\]
We obtain 
\begin{align*}
|hA| & = \left|  \bigcup_{i=0}^h L_i \right| 
 = \left| L_0^* \right|  + \sum_{i=i_0+2}^h  \left|  L_i \right| \\ 
  & =   (i_0+1)(b-a+c)  + (h-i_0)(ha + 1) \\
  & \qquad - \frac{(h+i_0+2)(h-i_0-1)(a-c)}{2} \\
  & = (i_0+1) b  +  (h-i_0)( ha + 1) - \frac{(h+i_0+1)(h-i_0)(a-c)}{2}. 
\end{align*}
This completes the proof. 
\end{proof}

\bc          \label{explicit:corollary:ab}
Let $h \geq 2$ and  $k \geq 3$.  

If $b> h(k-2)$, then the set $A= [0,k-2] \cup \{b\}$ satisfies
\[
|hA| =  (h+1)\left( 1 + \frac{h(k-2)}{2}  \right).
\] 

If $b \in [k-1, h(k-2)]$, then there exist unique integers 
 $i_0 \in [0,h-2]$ and $r \in [0,k-3]$ such that 
\[
b = (h-i_0)(k-2)-r 
\]
 and  the set $A= [0,k-2] \cup \{b\}$ satisfies 
\[ 
|hA| = (i_0+1) b +  (h-i_0)( h(k-2) + 1) - \frac{(h+i_0+1)(h-i_0)(k-2)}{2}.   
\] 
\ec

\begin{proof} 
This follows directly from Theorem~\ref{explicit:theorem:abc} with $a=k-2$ and $c=0$.  
\end{proof}

\bc           \label{explicit:corollary:ab2}
For all $h \geq 2$, 
\[
  (h+1)^2 \in \mcr_{\Z}(h,4) 
\] 
and, for all  
 $i_0 \in [0,h-2]$ and $r \in [0,1]$, 
\[
 (i_0+1) (2(h-i_0) -r ) +  (h-i_0)^2 \in \mcr_{\Z}(h,4).  
\] 
\ec

\begin{proof}
Apply Corollary~\ref{explicit:corollary:ab} with $k=4$.  
Let $A= [0,2] \cup \{b\}$.  
If $b> 2h$, then 
\[
|hA| =  (h+1)^2 \in \mcr_{\Z}(h,4).
\] 

If $i_0 \in [0,h-2]$ and $r \in [0,1]$, then 
\[
b = 2(h-i_0) -r \in  [3, 2h].  
\]
Conversely, if $b \in [3, 2h]$, then there exist unique integers 
 $i_0 \in [0,h-2]$ and $r \in [0,1]$ such that 
\[
b = 2(h-i_0) -r . 
\]
Setting $A = [0,2] \cup \{b\}$, we obtain   
\[ 
|hA| = (i_0+1) b +  (h-i_0)^2  \in \mcr_{\Z}(h,4). 
\] 
 This completes the proof. 
\end{proof}

\bc                     \label{explicit:corollary:ak}
For all $h \geq 2$ and  $k \geq 3$,   
\[
hk \in \mcr_{\Z}(h,k).
\]
\ec

\begin{proof}
Applying Corollary~\ref{explicit:corollary:ab} with $i_0 = h-2$ and $r = k-4$, 
we obtain $b=k$, $A = [0,k-2]\cup \{k\}$, and 
$|hA| = hk$. This completes the proof. 
\end{proof}

\bc                     \label{explicit:corollary:AP}
Let $h \geq 2$ and  $k \geq 3$.  
For all $i_0 \in [0,h-2]$, the sumset size set $\mcr_{\Z}(h,k)$  contains the arithmetic progression
\[
(i_0+1) b +  (h-i_0)( h(k-2) + 1) - \frac{(h+i_0+1)(h-i_0)(k-2)}{2}
\]
for $b \in [ (h-i_0)(k-2) - (k-3), (h-i_0)(k-2)]$.

In particular, $\mcr_{\Z}(h,k)$ contains the integer interval 
\[
\left[  \frac{ h^2 (k - 2)}{2} + \frac{hk}{2} - k + 3 ,    \frac{ h^2 (k - 2)}{2} + \frac{hk}{2}  \right].
\]
\ec

\begin{proof} 
This follows directly from Corollary~\ref{explicit:corollary:ab}.  
\end{proof}

For related work on sumset sizes, see~\cite{fox-krav-zhan25, hegy96,krav25,nath25aa,obry25,peri-roto25}.

\end{document}